\documentclass[10pt]{article}

\usepackage[margin=1in]{geometry}
\usepackage{amsmath}
\usepackage{cite}
\usepackage{mathtools}
\usepackage{amssymb}
\usepackage{amsthm}
\usepackage{graphicx}
\usepackage{color}
\usepackage[dvipsnames]{xcolor}
\usepackage{hyperref}
\usepackage{algorithm}
\usepackage{algpseudocode}
\usepackage[normalem]{ulem}
\usepackage{braket}
\usepackage{authblk}

\newcommand{\CC}{\mathbb{C}}

\newcommand{\OO}[1]{\mathcal{O}\paren{#1}}
\newcommand{\paren}[1]{\left( #1 \right)}
\newcommand{\brak}[1]{\left[ #1 \right]}

\newcommand{\wh}[1]{\widehat{#1}}
\newcommand{\wt}[1]{\widetilde{#1}}
\newcommand{\wb}[1]{\overline{#1}}

\renewcommand{\Im}{\operatorname{Im}}

\DeclareMathOperator\tr{Tr}

\newcommand{\GM}{G^M}
\newcommand{\GR}{G^R}
\newcommand{\GTV}{G^\rceil}

\newcommand{\SM}{\Sigma^M}
\newcommand{\SR}{\Sigma^R}
\newcommand{\STV}{\Sigma^\rceil}

\newcommand{\QTV}{Q^\rceil}
\newcommand{\gM}{g^M}
\newcommand{\qTV}{q^\rceil}
\newcommand{\sR}{\sigma^R}
\newcommand{\sTV}{\sigma^\rceil}
\newcommand{\gTV}{g^\rceil}
\newcommand{\muab}{\mu^{\text{AB}}}
\newcommand{\muam}{\mu^{\text{AM}}}
\newcommand{\mug}{\mu^{\text{G}}}

\newcommand{\Schrod}{Schr{\"o}dinger }

\providecommand{\keywords}[1]
{
  {\small \textbf{\textit{Keywords ---}} #1}
}

\title{A fast time domain solver for the equilibrium Dyson equation}

\author[1,2]{Jason Kaye \thanks{jkaye@flatironinstitute.org}}
\author[3]{Hugo U. R. Strand \thanks{hugo.strand@oru.se}}

\affil[1]{{\footnotesize Center for Computational Mathematics, Flatiron Institute, New York, NY 10010, USA}}
\affil[2]{{\footnotesize Center for Computational Quantum Physics, Flatiron Institute, New York, NY 10010, USA}}
\affil[3]{{\footnotesize School of Science and Technology, \"Orebro University, Fakultetsgatan 1, SE-701 82, \"Orebro, Sweden}}

\date{}

\begin{document}

\maketitle

\begin{abstract}

We consider the numerical solution of the real time equilibrium Dyson equation,
  which is used in calculations of the dynamical properties of quantum
  many-body systems. We show that
  this equation can be written as a system of coupled,
  nonlinear, convolutional Volterra integro-differential equations, for
  which the kernel depends self-consistently on the solution. As is
  typical in the numerical solution of Volterra-type equations, the
  computational bottleneck is the quadratic-scaling cost of history
  integration. However, the structure of the nonlinear Volterra integral
  operator precludes the use of standard fast algorithms. We propose a
  quasilinear-scaling FFT-based algorithm which respects the structure
  of the nonlinear integral operator. The resulting method can reach
  large propagation times, and is thus
  well-suited to explore quantum many-body phenomena at low energy
  scales. We demonstrate the solver with two standard model systems: the
  Bethe graph, and the Sachdev-Ye-Kitaev model.

\end{abstract}

\smallskip

\keywords{nonlinear Volterra integral equations, fast algorithms,
equilibrium Dyson equation, many-body Green's function methods}

\section{Introduction}

Quantum many-body systems \cite{Abrikosov:1975aa, Fetter:2003aa,
Negele:1998aa} can be described in terms of many-body Green's functions,
which characterize the system's response to the addition and removal of one or more
particles at different points in time.
The number of degrees of freedom required to describe a many-body
Green's function scales polynomially with the system size, independently of the number of particles.
Green's function
methods therefore provide an important complement to wavefunction-based
methods, especially for macroscopic systems, since the dimensionality of
a wavefunction is proportional to the number of particles it describes.
A large number of physical properties can be determined just from
single particle Green's functions, which give the expectation value
associated with the addition and removal of a single particle.

The equation of motion for the single particle Green's function is
called the Dyson equation \cite{Abrikosov:1975aa}.
The static thermal equilibrium properties of quantum many-body systems
can be obtained by solving the equation of motion in imaginary time
\cite{Negele:1998aa}, and the dynamical properties by solving it in
real time or frequency.
The last decade has seen a revival of time domain Green's function
methods \cite{Stefanucci:2013oq},
mainly driven by an interest in non-equilibrium phenomena. This has in turn spurred
advances in the development of numerical algorithms both for imaginary time
\cite{Dahlen:2005aa,boehnke11,gull18,dong20,ku00,ku02,kananenka16, shinaoka17,chikano18,li20,shinaoka21_2,kaye21_3}
and non-equilibrum real time \cite{Stan_2006, Stan:2009ab,
Aoki:2014kx,talarico19,Schuler:2020uy, kaye21_2} Green's functions.
However, the real time Dyson equation can also be used to determine the dynamical properties
of quantum many-body systems in equilibrium, by evolving the imaginary
time Green's function along the real time axis. 

As we will show in Section \ref{sec:dyson}, the Dyson equation of motion for the equilibrium
single particle Green's function in real time can be written as a system
of coupled nonlinear Volterra integro-differential equations of the form
\begin{equation}\label{eq:dysongeneral}
  \begin{cases}
    &i y'(t) + \int_0^t k\paren{y(t-t'),t-t'} y(t') \, dt' =
      f(y(t),t) \\
    &y(0) = y_0.
  \end{cases}
\end{equation}
In that context, the solution $y$ is a real time single particle
Green's function, and the interaction kernel $k$ is called the
self-energy.
For simplicity of exposition, we take all quantities to be complex scalar-valued, but the matrix-valued case often encountered in Green's function calculations is a straightforward generalization.
We also note that although we have assumed
local nonlinearities $k(t) = k(y(t),t)$ and $f(t) = f(y(t),t)$, our method
applies equally well to the more general case of causal nonlinearities
$k(t) = k(y(t')|_{t' \leq t},t)$ and $f(t) = f(y(t')|_{t' \leq t},t)$, which appear in
more complicated diagrammatic models of the self-energy.

The cost of computing the history integral term in
\eqref{eq:dysongeneral} scales as $\OO{N^2}$ with the number $N$ of time
steps, severely limiting the propagation times achievable by time domain
solvers \cite{Strand:2015ac}.
This creates a practical barrier in studies of collective low energy
phenomena in quantum many-body systems, which require high resolution in
the frequency domain.
This approach is typically avoided as a result.
The standard alternative is to evaluate the self-energy in the time
domain, where it is typically simpler, and transform to the frequency
domain to solve the Dyson equation and obtain the corresponding Green's function.
While this method eliminates the problem of history integration, and can be
accelerated using the fast Fourier transform (FFT) \cite{Foerster:2011ty}, it requires solving
a global nonlinear problem
\cite{PhysRevB.54.8411,PhysRevB.57.2108}. The convergence properties of
a given nonlinear iteration procedure are problem dependent, and the
number of iterations can grow in physically interesting regimes, such as
the low temperature limit.

We also mention another method, which is to perform analytic continuation
to obtain a simpler integral equation in the imaginary time domain, and then
analytically continue the solution back to the real time or frequency
domain by solving a severely ill-conditioned integral equation \cite{anacont}.
For sufficiently complicated approximations to the self-energy, direct solution
in real time or frequency can be impractical, and this technique is the only
option. However, despite significant recent progress in sophisticated
regularization methods \cite{anacontNevanlinna, anacontCaratheodory, huang23}, the
analytic continuation problem is fundamentally ill-posed, and obtaining high
quantitative accuracy reliably, as we seek here, remains a challenge. 

In this work, we present a time stepping algorithm for \eqref{eq:dysongeneral}
which reduces the cost of Volterra history integration to $\OO{N \log^2
N}$, eliminating the main barrier to carrying out equilibrium
Green's function calculations in the time domain. Using this approach, the
nonlinear problem is local and well-controlled; the number of
nonlinear iterations is independent of global features of the solution,
and is typically very small.
In combination with the recently introduced discrete Lehmann
representation \cite{kaye21_3} of the imaginary time axis, our time
domain algorithm enables the calculation of dynamical properties of
quantum many-body systems at low temperatures and large
propagation times, i.e.\ at low energy scales, which are
essential in capturing emergent many-body excitations.

The problem of efficient Volterra history integration is well known in
the scientific computing literature, and
several techniques have been proposed, particularly
in the context of solving Volterra integral
equations and computing nonlocal transparent boundary conditions. We
mention fast Fourier transform (FFT)-based methods \cite{hairer85},
methods based on sum-of-exponentials projections of the history
\cite{alpert00,jiang04,jiang08,jiang15,veerapaneni07,li10,wang19,kaye21,kaye21_4}, convolution
quadrature \cite{lubich02,schadle06}, and hierarchical low-rank or
butterfly matrix compression \cite{kaye21_2,dolz21,kaye18}.
Of the methods mentioned above, several have been applied to the
efficient numerical
solution of convolutional nonlinear Volterra integral equations
\cite{hairer85,schadle06,dolz21}.
However, we emphasize that the form of the nonlinearity in \eqref{eq:dysongeneral}
is different from that typically considered,
\begin{equation} \label{eq:nlvio}
  \int_0^t k(t-t') F\paren{y(t'),t'} \, dt',
\end{equation}
where $k$, or its Laplace transform, is given either explicitly or
numerically. All of the methods cited above make use of this \textit{a
priori}
access to $k$ in order to obtain a fast history summation scheme, and are
therefore not applicable in the present setting, in which $k$ itself is
determined during the course of time stepping by a self-consistency
condition. This can be made evident by considering a simple example, $k(y(t),t) = y(t)$, which yields the
integral operator $\int_0^t y(t-t') y(t') \, dt'$ in
\eqref{eq:dysongeneral}. We refer to the form of the nonlinearity
considered here as a kernel nonlinearity, to distinguish it from
nonlinearities of the form \eqref{eq:nlvio}.

Our algorithm is inspired by the FFT-based method of Hairer,
Lubich, and Schlichte (HLS), proposed for the numerical solution of Volterra
integral equations with nonlinearities of the typical form \eqref{eq:nlvio}
\cite{hairer85}. Their approach, which will be discussed in detail in
Section \ref{sec:hls}, is to restructure the history
sums into a hierarchy of Toeplitz block matrix-vector products, and to apply the blocks in
quasi-linear time using the FFT. This yields an $\OO{N \log^2 N}$
algorithm. Although their method does not
apply directly to the present case because of the kernel nonlinearity,
we find that a more sophisticated block partitioning gives a compatible
algorithm with the same computational complexity.

This paper is organized as follows. We begin in Section \ref{sec:fastalg} by
describing our fast algorithm to compute the history integrals in
\eqref{eq:dysongeneral}. In Section \ref{sec:mbt}, we give a brief
introduction to many-body Green's functions and the equilibrium Dyson
equation for single particle Green's functions, and show that the latter
can be reduced to a system of coupled equations of the form
\eqref{eq:dysongeneral}.
In Section \ref{sec:examples}, we apply our method to two canonical model systems:
the Bethe graph \cite{QPGFBook, Georges:1996aa} and the
Sachdev-Yitaev-Ke (SYK) model \cite{PhysRevLett.70.3339,chowdhury21,Gu:2020vl}.
We demonstrate a
significant increase in attainable propagation times, and therefore
improvements in the resolution of the fine spectral properties of
solutions. A concluding discussion is given in Section
\ref{sec:conclusion}. We also
provide appendices describing a specific high-order time stepping method
for \eqref{eq:dysongeneral}, and technical algorithmic details of
the proposed method.

\section{Fast history summation} \label{sec:fastalg}

To illustrate the appearance of the history sums in a simple manner, we
first write down a forward Euler discretization of
\eqref{eq:dysongeneral}. We do not recommend using this method in
practice, due to its poor accuracy and stability properties. Rather, we
propose a high-order implicit multistep method in Appendix
\ref{app:adams}. History sums of the same form appear no matter the
choice of discretization, and the details of our fast history summation
algorithm are independent of this choice.

Applying the forward Euler method to \eqref{eq:dysongeneral}
yields
\[y_{n+1} = y_n + i \Delta t \int_0^{t_{n}} k(t_{n}-t')
y(t') \, dt' - i \Delta t f_{n}.\]
Here $t_n = n \Delta t$ for a time step $\Delta t$ and $n = 0,\ldots,N-1$,
$y_n \approx y(t_n)$, and $f_n = f(t_n)$. For notational simplicity, we have suppressed
the dependence of $k$ and $f$ on $y$. Applying the trapezoidal rule to the integral
gives
\[y_{n+1} = y_n + i \Delta t^2 \sum_{m=0}^n
k_{n-m} y_m -\frac{i \Delta t^2}{2} \paren{k_n y_0 + k_0 y_n}  - i \Delta t f_{n}.\]
The computational bottleneck is evidently the calculation of the history
sum
\begin{equation} \label{eq:histsum}
  s_n = \sum_{m=0}^{n} k_{n-m} y_m
\end{equation}
at each time step. Indeed, the cost to compute all history sums for $n =
0,\ldots,N-1$ is $\OO{N^2}$,
and dominates the otherwise $\OO{N}$ cost of time stepping.

The collection of history sums takes the form of a
lower triangular Toeplitz matrix-vector product $s = Ky$. Here, $s, y \in
\CC^N$ and $K \in \CC^{N \times N}$, with $K_{nm} = k_{n-m}$ for $n \geq
m$ and $K_{nm} = 0$ for $n < m$, $n,m = 0,\ldots,N-1$.
An $N \times N$ Toeplitz matrix can be
applied to a vector in $\OO{N \log N}$ operations using the FFT;
a detailed description of the algorithm is given in
Appendix \ref{app:toeplitz}.
Therefore, if $y$ and $k$
were known, then computing the sums \eqref{eq:histsum} with quasi-linear cost would
be straightforward.

However, in the present setting, there are two additional complications:
\begin{enumerate}
  \item $y_n$ must be computed in order to obtain $s_n$, which must in
    turn be computed in order to obtain $y_{n+1}$.
  \item $y_n$ must be computed in order to obtain $k_n$, which must in
    turn be computed in order to obtain $y_{n+1}$.
\end{enumerate}

The first, which is familiar in the literature on fast algorithms for
Volterra-type equations and integral operators, prevents us from
computing the history sums all at once as a matrix-vector product, so we
cannot simply use the standard algorithm for the fast application of
Toeplitz matrices. The HLS algorithm \cite{hairer85}, which we review in
Section \ref{sec:hls}, offers a solution to this problem, as do the many
other methods for the efficient evaluation of Volterra integral
operators mentioned in the introduction. The second problem has not, to
our knowledge, been addressed in this literature. We will propose a solution in Section
\ref{sec:fullalg}, which appropriately modifies the HLS algorithm to
respect the causal structure of the kernel nonlinearity.

\subsection{HLS algorithm for fast history summation with a fixed kernel}
\label{sec:hls}

\begin{figure}[t]
  \centering
  \includegraphics[width=0.45\textwidth]{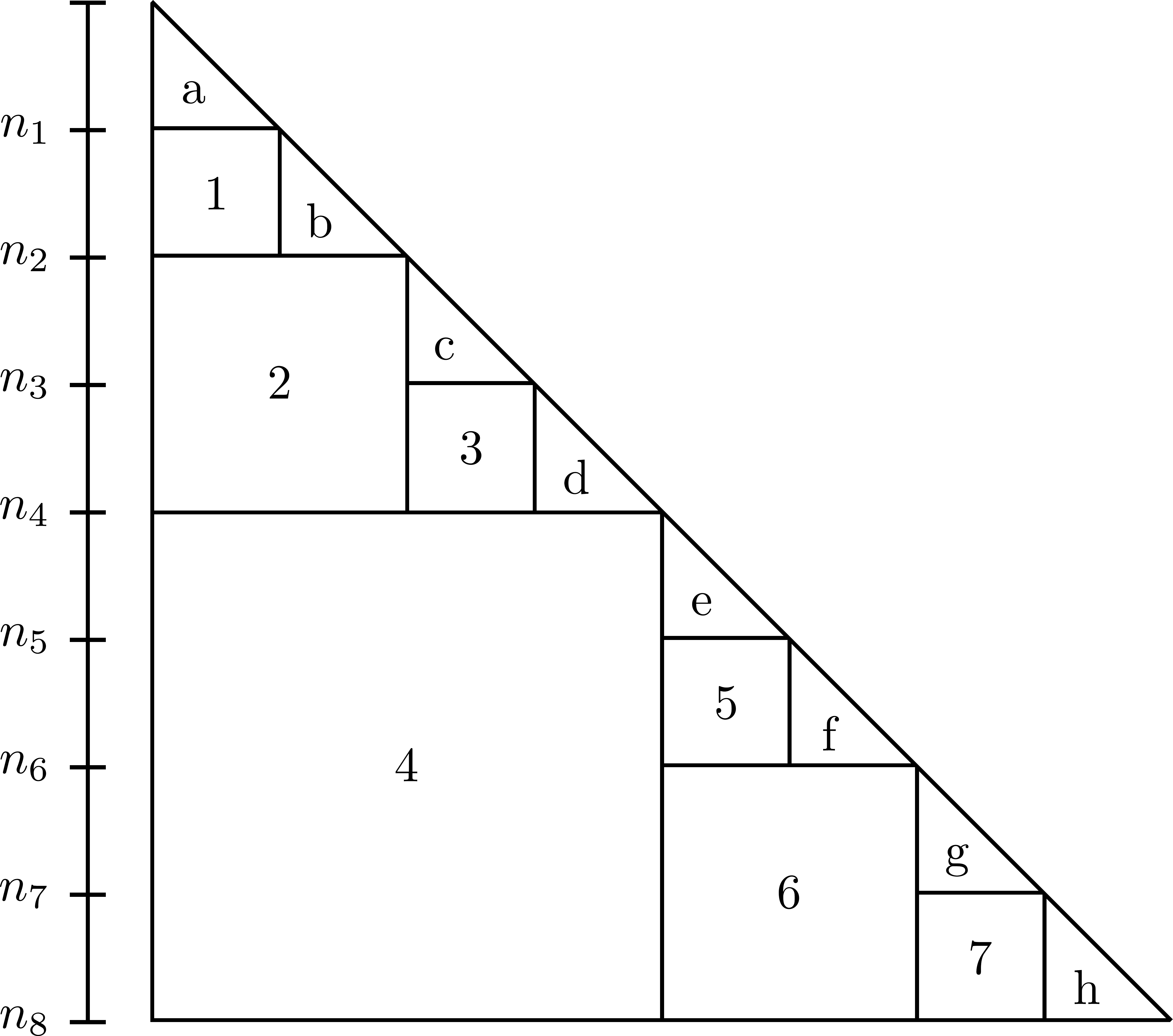}
  \caption{Illustration of the HLS algorithm to compute the history
  sums \eqref{eq:histsum} on the fly with a given kernel $k$. Here, we
  have $L = 3$ levels. The
  lettered blocks correspond to entries which are incorporated into the
  partial sums by direct summation, and the numbered blocks are applied
  using a fast Toeplitz matrix-vector product. The blocks are processed
  in the order corresponding to their first row indices: a-1-b-2-c-3-d-4-e-5-f-6-g-7-h.}
  \label{fig:triblock}
\end{figure}

The HLS algorithm computes the sums
\eqref{eq:histsum}, for the case in which the kernel $k$ is known, in
$\OO{N \log^2 N}$ operations \cite{hairer85}. The matrix
$K_{nm} = k_{n-m}$ is partitioned into square blocks as in Figure
\ref{fig:triblock}.
In particular, we define $n_j = \lceil jN/2^L \rceil$ for $j =
1,\ldots,2^L-1$ and $n_{2^L} = N-1$, where $L$ is the number of levels
of subdivision. We choose $L = \OO{\log_2 N}$ so that $n_1$ is a small,
fixed constant. Each block is Toeplitz, and
therefore can be applied with quasi-optimal complexity using the FFT.

The algorithm proceeds by accumulating partial sums $\wt{s}_n$ in
an efficient manner, such that by the $n$th time step, $\wt{s}_n =
s_n$ as needed. Contributions to the partial sums come both from
triangular blocks (shown as the lettered blocks in
Figure \ref{fig:triblock}) applied directly row-by-row, and square blocks (shown as the numbered
blocks in Figure \ref{fig:triblock}) applied using the FFT. The
algorithm amounts to a re-ordering of the calculation of history
sums from the standard row-by-row approach, so as to enable the use of the
FFT.

A pseudocode is given in Algorithm \ref{alg:hls}. We note the use of
colon notation for indexing into vectors; for a vector $x$, we define
$x_{i:j} \equiv (x_i,\ldots,x_j)^T$. Let us go through the
first few steps of the algorithm.
\begin{algorithm}
  \begin{algorithmic}[1]
    \State \textbf{Given:} The first and last row indices, $r_j^f$ and
      $r_j^l$, and the first and last column indices, $c_j^f$ and
      $c_j^l$, of block $j$ in Figure \ref{fig:triblock}, for $j =
      1,\ldots,2^L-1$
    \State Initialize $\wt{s}_{0:N-1} = 0$
    \For{$n=0,n_1-1$}
      \State $s_n \gets \sum_{m=0}^n k_{n-m} y_m$
      \State Take a time step to obtain $y_{n+1}$
    \EndFor
    \For{$j=1,2^L-1$}
      \State Apply block $j$ to $y_{c_j^f:c_j^l}$ using FFT-based
      algorithm, and add result to $\wt{s}_{r_j^f:r_j^l}$
      \For{$n=n_j,n_{j+1}-1$}
        \State Compute local contribution to history sums: $\wt{s}_n \gets \wt{s}_n + \sum_{m=n_j+1}^n k_{n-m} y_m$
        \State Set $s_n = \wt{s}_n$ and take a time step to obtain
        $y_{n+1}$
      \EndFor
    \EndFor
  \end{algorithmic}
  \caption{Time-stepping using the HLS algorithm}
  \label{alg:hls}
\end{algorithm}

First, the time steps $n = 0,\ldots,n_1-1$ are carried out by the standard
method, with the history sums $s_0,\ldots,s_{n_1-1}$ computed by direct
summation. Once $y_{n_1}$ has been computed,
block 1 in Figure \ref{fig:triblock} can be applied to the vector
$y_{0:n_1}$, yielding partial contributions to the sums
$s_{n_1:n_2-1}$. These partial contributions are stored in
$\wt{s}_{n_1:n_2-1}$.

Next, we carry out time steps $n = n_1$ to
$n = n_2-1$. At the $n$th time step, we add the local contribution corresponding
to the $n$th row of block b to the stored partial sum $\wt{s_n}$:
$\wt{s}_n \gets \wt{s}_n + \sum_{m=n_1+1}^n k_{n-m} y_m$. We then take
$s_n = \wt{s}_n$, and can compute $y_{n+1}$.

Once $y_{n_2}$ has been obtained, block 2 can be applied to $y_{0:n_2}$,
and the result stored in $\wt{s}_{n_2:n_4-1}$. Then we proceed with time
steps $n = n_2$ to $n = n_3-1$. At the $n$th time step, we add the local
contribution corresponding to the $n$th row of block c to $\wt{s}_n$,
$\wt{s}_n \gets \wt{s}_n + \sum_{m=n_2+1}^n k_{n-m} y_m$, set $s_n =
\wt{s}_n$, and obtain $y_{n+1}$.

Once $y_{n_3}$ is obtained, block 3 can be applied to $y_{n_2+1:n_3}$,
and its contribution added to $\wt{s}_{n_3:n_4-1}$. Then the local
contributions corresponding to the rows of block d are added to the
partial sums in order to take the corresponding time steps.
The algorithm proceeds in this manner, with contributions from each
block added to the corresponding partial sums in the indicated order,
such that by the end of the $n$th time step, $\wt{s_n} = s_n$ as needed.

To derive the computational complexity of this algorithm, we simply add
up the cost of applying all blocks. Recall that $L = \OO{log_2 N}$. Since each triangular block is of
dimensional approximately $N/2^L = \OO{1}$, and there are $2^L = \OO{N}$
such blocks, the total cost of applying these blocks is $\OO{N}$. For
the square blocks, we have one block of dimension approximately $N/2$,
two blocks of dimension $N/4$, four blocks of dimension $N/8$, and so
on. The total cost of applying these blocks using the fast Toeplitz
algorithm is therefore approximately
\[\sum_{l=1}^L 2^{l-1} \times \frac{N}{2^l} \log\paren{\frac{N}{2^l}} =
\OO{N \log^2 N}.\]

\subsection{Fast history summation with a kernel nonlinearity} \label{sec:fullalg}

The HLS algorithm cannot be applied
if $k_{n+1},\ldots,k_{N-1}$ are unknown at the
$n$th time step. Indeed, the block $j$ is applied as soon as $y_{n_j}$ is
obtained. However, block $j$ may contain values $k_n$ of the kernel with
$k > n_j$; for example, block $4$ contains values up to $k_{N-1}$. Since
$k_n$ depends on $y_n$, it is not possible to apply the full blocks in
the order prescribed by the algorithm.
In particular, only the upper triangular parts of blocks 1, 2, and 4
are known at the step they must be applied.

A modified partitioning of the same matrix is depicted in Figure
\ref{fig:triblockmod}. It contains square, triangular, and
parallelogram-shaped blocks. In Appendix \ref{app:toeplitz}, we generalize the standard
FFT-based algorithm for Toeplitz matrices to triangular and
parallelogram-shaped blocks (appropriately completed to full rectangular matrices
by zero-padding), allowing these blocks to be applied in $\OO{n \log n}$
operations as well.

\begin{figure}[t]
  \centering
  \includegraphics[width=0.45\textwidth]{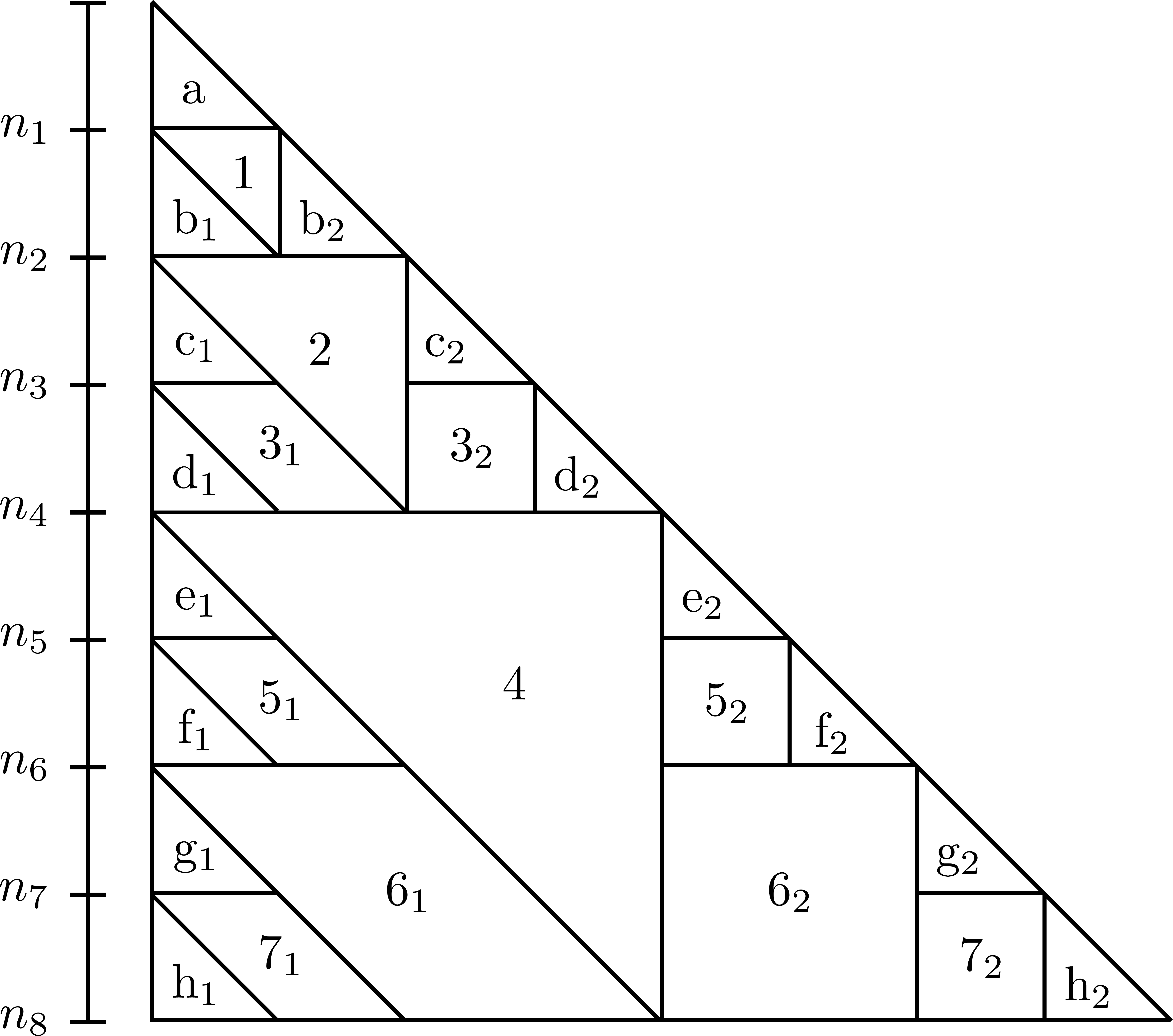}
  \caption{Modification of the block structure depicted in Figure
  \ref{fig:triblock} to account for kernels which are obtained
  on the fly. The blocks are again processed
  in the order corresponding to their first row indices:
  a-1-b-2-c-3-d-4-e-5-f-6-g-7-h. Blocks with the same label and
  subscripts $1$ and $2$, like $3_1$ and $3_2$, can be processed simultaneously.}
  \label{fig:triblockmod}
\end{figure}

\begin{figure}[t]
  \centering
  \includegraphics[width=0.65\textwidth]{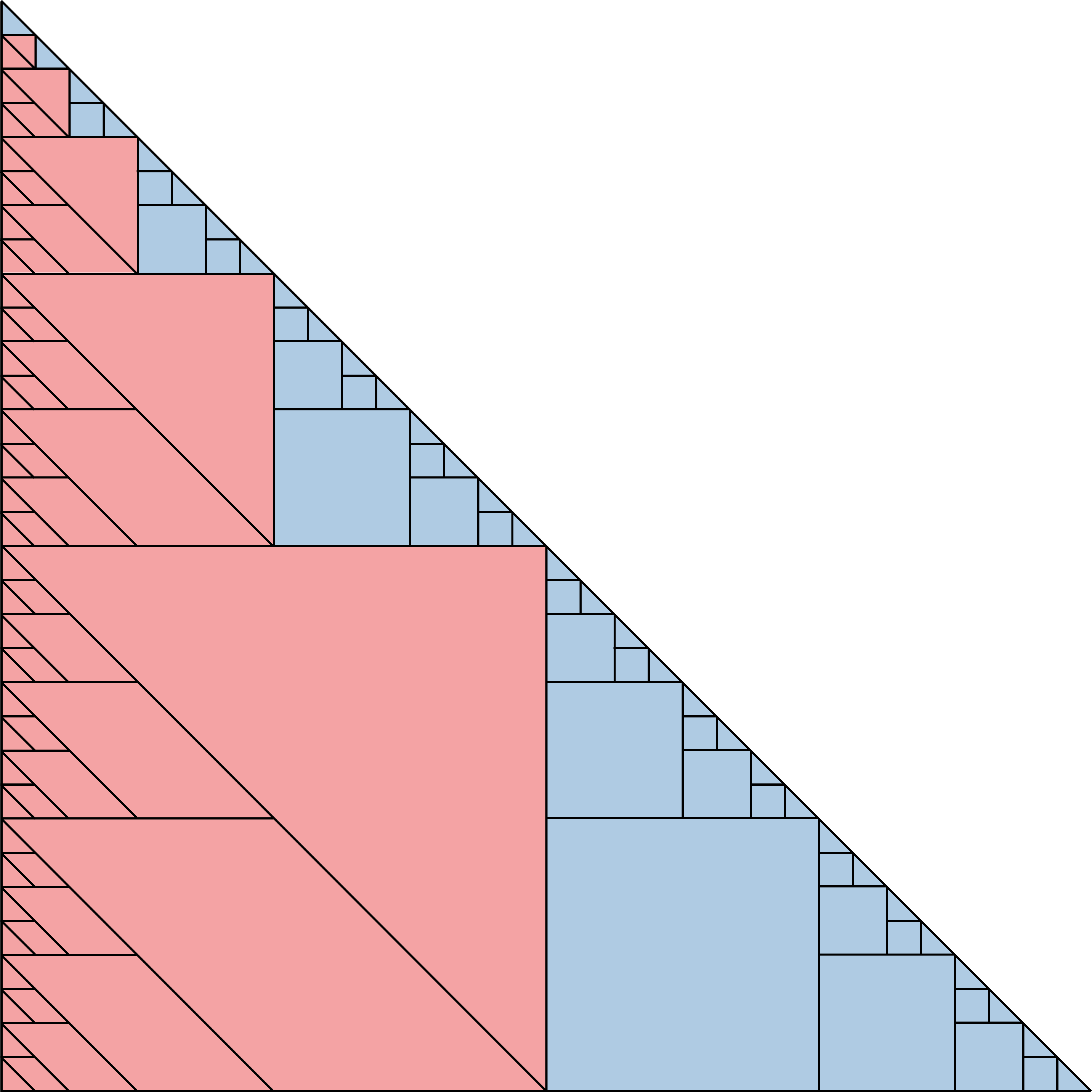}
  \caption{Illustration of the block partitioning of the kernel matrix
  for the full fast algorithm, with $L = 5$. The full partioning is
  comprised of triangular units (blue) from the HLS algorithm (Figure
  \ref{fig:triblock}), and square units (red), which account for the
  self-consistent nature of \eqref{eq:dysongeneral}.}
  \label{fig:tribig}
\end{figure}

The structure of our algorithm is similar to that of the HLS algorithm,
but with the modified partitioning. The triangular blocks indicated in
Figure \ref{fig:triblockmod} by letters contain rows which are applied
directly, and the numbered blocks are applied as soon as the time step
corresponding to their first row
$n_j$ is reached. The shapes of the numbered blocks ensure that their
entries depend only on $k_n$ for $n \leq n_j$. The results of both the
direct and block applications are accumulated in partial sums $\wt{s}_n$
as before. Unlike in the HLS algorithm, at each time step there are two
local contributions that must be made; the first, as before,
corresponds to values $y_m$ for $m$ near the current time step $n$,
and the second to values $k_m$ for $m$ near $n$. Also, at some steps,
two blocks are applied instead of one. A pseudocode for the full
procedure is given in Algorithm \ref{alg:hlsmod}.
The partioning for $L = 5$ is shown in Figure \ref{fig:tribig}. This is
similar to Figure \ref{fig:triblockmod}, but contains more levels.

\begin{algorithm}
  \begin{algorithmic}[1]
    \State \textbf{Given:} The first and last row indices, $r_j^f$ and
      $r_j^l$, and the first and last column indices, $c_j^f$ and
      $c_j^l$, of block $j$ in Figure \ref{fig:triblock}, for $j =
      1,2,3_1,3_2,4,5_1,5_2,\ldots$
    \State Initialize $\wt{s}_{0:N-1} = 0$
    \For{$n=0,n_1-1$}
      \State $s_n \gets \sum_{m=0}^n k_{n-m} y_m$
      \State Take a time step to obtain $y_{n+1}$
    \EndFor
    \For{$j=1,2^L-1$}
      \If{there is one numbered block at row $n_j$}
        \State Apply block $j$ to $y_{c_j^f:c_j^l}$ using FFT-based
          algorithm, and add result to $\wt{s}_{r_j^f:r_j^l}$
      \Else
        \State There are two blocks at row $n_j$; do the same as above, for $j = j_1$ and $j = j_2$
      \EndIf
      \For{$n=n_j,n_{j+1}-1$}
        \State Compute local contribution to history sums: $\wt{s}_n \gets
        \wt{s}_n + \sum_{m=0}^{n-n_j-1} k_{n-m} y_m + \sum_{m=n_j+1}^n k_{n-m} y_m$
        \State Set $s_n = \wt{s}_n$ and take a time step to obtain
        $y_{n+1}$
      \EndFor
    \EndFor
  \end{algorithmic}
  \caption{Time-stepping using the fast
  history summation algorithm for kernel nonlinearities}
  \label{alg:hlsmod}
\end{algorithm}

To obtain the computational complexity of the new algorithm, we consider
the two types of units separately. We know from the previous section
that the total cost of applying the triangular HLS-type units,
indicated in blue in Figure \ref{fig:tribig}, is approximately
\[\sum_{l=1}^{L} \frac{N}{2^l} \log^2\paren{\frac{N}{2^l}} = \OO{N \log^2 N}.\]
The total cost of applying each block in the square unit of dimension $N/2^l$
is
\[\frac{N}{2^l} \log\paren{\frac{N}{2^l}} + \sum_{k=0}^{L-l} 2^k
  \frac{N}{2^{l+k}} \log\paren{\frac{N}{2^{l+k}}} \\ \leq \frac{N}{2^l}
  \log\paren{\frac{N}{2^l}} (L-l+2).\]
Here, the first term corresponds to the top-right upper triangular
block, and the second term corresponds to the rest of the blocks.
Summing over all square units gives the total cost
\[ \sum_{l=1}^L \paren{\frac{N}{2^l}
  \log\paren{\frac{N}{2^l}} (L-l+2)} = \OO{N \log^2 N}.\]
Thus the computational complexity of the modified algorithm is $\OO{N
\log^2 N}$, the same as that of the original HLS algorithm.

\section{Green's functions and the equilibrium Dyson equation} \label{sec:mbt}

The aim of this section is to give a brief introduction to quantum
many-body theory, the single particle Green's function, and the Dyson
equation of motion, starting from basic quantum mechanics. We refer the
reader to \cite{Stefanucci:2013oq} for a more thorough introduction to
the subject.

In thermal equilibrium, the state of a quantum many-body
system is described by a many-body density matrix $\hat{\varrho}$, which
represents a statistical mixture of states $\ket{\psi_n}$ weighted,
according to the Gibbs distribution, by $e^{-\beta E_n}/Z$:
\begin{equation*}
  \hat{\varrho} = \frac{e^{-\beta \hat{H}}}{Z} = \frac{1}{Z} \sum_n
  e^{-\beta E_n} | \psi_n \rangle \langle \psi_n |.
\end{equation*}
Here $\hat{H}$ is the Hamiltonian of the system for times $t < 0$, $\beta$ is the inverse temperature, $|\psi_n\rangle$ is an eigenstate of the system with energy $E_n$, $\hat{H} |\psi_n \rangle = E_n | \psi_n \rangle$, and $Z$ is the partition function $Z = \sum_n e^{-\beta E_n}$.

In this context, the expectation value of an operator $\hat{O}$ is given by the ensemble average
\[
\braket{\hat{O}} =
\frac{1}{Z} \sum_n e^{-\beta E_n}  \bra{\psi_n} \hat{O} \ket{\psi_n}
= \text{Tr} \left[ \hat{\varrho} \, \hat{O} \right] .
\]

For a system with a time-dependent Hamiltonian $\hat{H}(t)$ for $t \ge 0$,
the operator expectation values are time-dependent, and can be
calculated in the Heisenberg picture by time evolving the operator: $\hat{O} \mapsto \hat{U}(0, t) \, \hat{O} \, \hat{U}(t, 0) \equiv \hat{O}(t)$. Here $\hat{U}(t,t')$ is the unitary time evolution operator from time
$t'$ to $t$ --
the solution operator for the \Schrod equation with Hamiltonian $\hat{H}(t)$
-- and is given formally by 
\[\hat{U}(t, t') =
  \begin{cases} 
    \mathcal{T} \exp\brak{-i \int_{t'}^t \hat{H}(\bar{t}) \, d\bar{t}} &
    \text{if } t > t'  \medskip  \\
    \wb{\mathcal{T}} \exp\brak{+ i \int_{t}^{t'} \hat{H}(\bar{t}) \, d\bar{t}} &
    \text{if } t < t',
  \end{cases}
\]
with $\mathcal{T}$ the time ordering operator and $\wb{\mathcal{T}}$ the
anti-time ordering operator.
The time-dependent ensemble average is then given by
\begin{equation} \label{eq:tdea}
  \braket{\hat{O}(t)} \equiv
  \tr\brak{\hat{\varrho} \, \hat{U}(0,t) \, \hat{O} \, \hat{U}(t,0)}.
\end{equation}
\begin{figure}[t]
  \centering
  \includegraphics[scale=1]{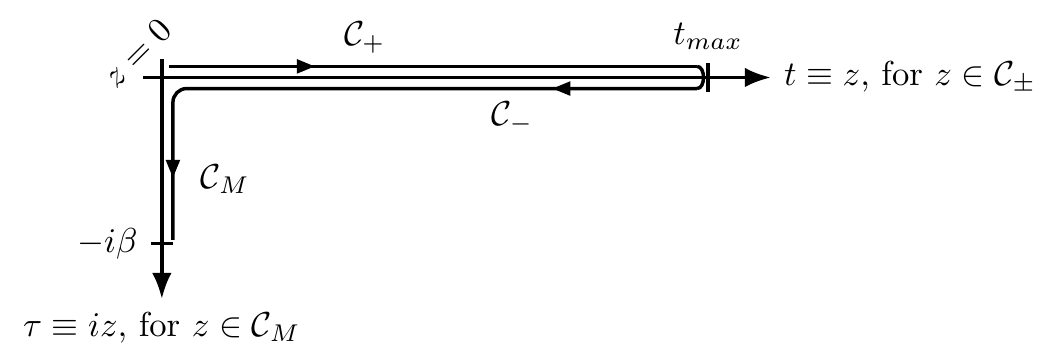}
  \caption{Generalized time contour for time-dependent correlation
  functions of quantum many-body systems. Here, $t_{max}$ indicates a finite
  propagation time which can in principle be taken to $\infty$.}
  \label{fig:contour}
\end{figure}
The notion of time propagation can also be extended to the initial
state, since $\hat{\varrho}$ can be expressed as
\begin{equation*}
  \hat{\varrho} = \frac{1}{Z} \exp \left[ - \int_0^{\beta} \hat{H} \,
  d\tau \right] = \frac{1}{Z} \exp \left[ -i \int_0^{-i\beta}
  \hat{H} \, dz \right] \equiv \frac{1}{Z} \hat{U}(-i\beta, 0).
\end{equation*}
Hence, the initial thermal equilibrium density matrix can be viewed in
terms of evolution in imaginary time from $t=z=0$ to $-i\beta$,
and the time-dependent ensemble average \eqref{eq:tdea} can be expressed as
\begin{equation}
  \braket{\hat{O}(t)} \equiv
  \frac{1}{Z} \tr\brak{ \hat{U}(-i\beta, 0) \, \hat{U}(0,t) \, \hat{O} \, \hat{U}(t,0)}
  .
  \label{eq:tdeaimtime}
\end{equation}
One can now interpret the contents of the trace in \eqref{eq:tdeaimtime} as the following
sequence of operations: i) propagate from time $0$ to
$t$, and apply the operator $\hat{O}$; ii) propagate
backwards from time $t$ to $0$; and iii) propagate in imaginary
time from $0$ to $-i\beta$.

More generally, expectation values of multiple operators acting at different
times, like the two-time correlation function between operators
$\hat{O}_1$ and $\hat{O}_2$,
\[
\langle \hat{O}_1(t) \hat{O}_2(t') \rangle
\equiv
\frac{1}{Z} \tr\brak{ \hat{U}(-i\beta, 0) \, \hat{U}(0,t) \, \hat{O}_1
\, \hat{U}(t,t') \, \hat{O}_2 \, \hat{U}(t',0)},
\]
can be treated by means of propagation on the generalized time contour $\mathcal{C}$
shown in Figure \ref{fig:contour}.
This contour is comprised of three branches, representing the three
directions of propagation used in calculating expectation values: a forward time propagation branch $\mathcal{C}_+$,
a backward time propagation branch $\mathcal{C}_-$, and an
imaginary time propagation branch $\mathcal{C}_M$.

\subsection{Single particle Green's functions} \label{sec:gf}

Quantum many-body systems can be described using a special class of
expectation values called many-body Green's functions. Many-body Green's
functions are the expectation values associated with the addition and
removal of particles from a system at a collection of times on
$\mathcal{C}$. For example, the single particle Green's function $G(z, z')$ is given by \cite{Stefanucci:2013oq}
\begin{equation} \label{eq:G}
  G(z, z') = -i \langle \mathcal{T}_\mathcal{C} c(z) c^\dagger(z') \rangle
  \, , \quad z, z' \in \mathcal{C}
  \, ,
\end{equation}
where $c^\dagger(z')$ is a particle creation operator, which adds a
particle to the system at a
generalized time $z'$, $c(z)$ is a particle annihilation operator,
which removes a particle from the system at another generalized time
$z$, and $\mathcal{T}_\mathcal{C}$ is the contour time ordering operator, which
commutes operators according to their contour ordering; see Figure \ref{fig:contour}.
In the standard treatment \cite{Stefanucci:2013oq}, Green's function components are introduced for each possible restriction of $z$ and $z'$ to the three branches $\mathcal{C}_+$, $\mathcal{C}_-$, and $\mathcal{C}_M$,
using the real-valued arguments $t$ and $\tau$, where $t \equiv z$ for
$z \in \mathcal{C}_\pm$, and $\tau \equiv iz$ for $z \in \mathcal{C}_M$;
see also Figure \ref{fig:contour}.

\subsection{The real time equilibrium Dyson equation} \label{sec:dyson}

The Dyson equation of motion for the single particle Green's function
\eqref{eq:G} can be derived from the Martin-Schwinger hierarchy of
coupled equations for n-particle Green's functions  by collecting
higher-order many-body corrections into an integral kernel called the
self-energy $\Sigma$. Here we write the Dyson equations
of motion for the Green's function components relevant in the case of
equilibrium real time evolution.
We note that in equilibrium, the Green's function is translation invariant in real time.

The initial thermal equilibrium state is described by the restriction
$G^M$ of the contour Green's function $G$ to imaginary time $z =
-i\tau$,
\begin{equation*}
  z, z' \in \mathcal{C}_M
  \, , \quad
  G(z, z') =
  G(-i\tau, -i\tau') =
  i G^M(\tau, \tau') =
  i G^M(\tau - \tau').
\end{equation*}
$\GM$ is referred to as the
Matsubara Green's function.
The Dyson equation of motion for $\GM(\tau)$, for $\tau \in [0,\beta]$, has the form
\begin{equation} \label{eq:dysonm}
  \begin{cases}
    &\paren{-\partial_\tau - h} \GM(\tau) - \int_0^\beta \SM(\tau-\tau')
  \GM(\tau') \, d\tau' = 0 \\
    &G^M(0) - \xi G^M(\beta) = -1.
  \end{cases}
\end{equation}
Here, $\xi = \pm 1$ for bosonic and fermionic particles, respectively,
and $\GM$ is extended to $(-\beta,0)$ in the convolution by the
periodicity or antiperiodicity conditions $\GM(-\tau) = \xi
\GM(\beta-\tau)$.\footnote{The boundary condition in \eqref{eq:dysonm}
is obtained from the periodicity/antiperiodicity condition and the commutation relation $G(0^-) - G(0^+) = 1$.}
The interaction kernel in these equations is the imaginary time self-energy $\Sigma^M$,
which depends self-consistently on $G^M$.

In equilibrium systems, it is sufficient to solve for the Green's function restricted to mixed real and imaginary time arguments.
The left mixing component $\GTV(t, \tau)$, with the second argument on
the vertical imaginary time branch $\mathcal{C}_M$,
\begin{equation*}
  z \in \mathcal{C}_{\pm} , \, z' \in \mathcal{C}_M
  \, , \quad
  G(z, z') =
  G(t, -i \tau) =
  \GTV(t, \tau)
  \, ,
\end{equation*}
satisfies the equilibrium Dyson equation of motion
\begin{equation} \label{eq:dysonlm}
  \begin{cases}
    &\paren{i\partial_t-h} \GTV(t,\tau) - \int_0^t \SR(t-t')
    \GTV(t',\tau) \, dt' = \QTV(t,\tau) \\
    &\GTV(0,\tau) = i \GM(-\tau) = i \xi \GM(\beta-\tau),
  \end{cases}
\end{equation}
with
\begin{equation} \label{eq:qlm}
  \QTV(t,\tau) = \int_0^\beta \STV(t,\tau') \GM(\tau'-\tau) \, d\tau'.
\end{equation}
Here, $\STV(t,\tau)$ and $\SR(t)$ are called the left mixing and
retarded self-energies, respectively, and they depend
self-consistently on $\GTV(t,\tau)$.
We note that the equation \eqref{eq:dysonm} for the Matsubara Green's
function can be solved first, independently, and its solution then used
in the initial condition and source term in \eqref{eq:dysonlm}.

The solution of \eqref{eq:dysonm} and \eqref{eq:dysonlm} completely
determines the single particle Green's function. However, several other
components are commonly used in the literature, and we mention them for
completeness. 
The lesser component is defined as
$G^<(t, t') = G(z, z')$ with $z \in \mathcal{C}_-$, $z' \in
\mathcal{C}_+$, and it is related to the occupied density of states. It
is given in equilibrium by
$G^<(t,t') \equiv G^<(t - t') = \GTV(t - t', 0)$.
The greater component is defined as
$G^>(t, t') = G(z, z')$ with $z \in \mathcal{C}_+$, $z' \in
\mathcal{C}_-$, and it is related to the unoccupied density of states.
It is given by
$G^>(t,t') \equiv G^>(t - t') = \xi \GTV(t - t', \beta)$.
The retarded component is given by
$G^R(t, t') = G^>(t, t') - G^<(t, t')$,
and is related to the full density of states. Again taking $\GR(t,t')
\equiv \GR(t-t')$, we therefore have
\begin{equation} \label{eq:grfromgtv}
  \GR(t) = -\paren{\GTV(t,0)-\xi\GTV(t,\beta)}.
\end{equation}
For later convenience, we note that $\GR(t)$ satisfies the Dyson equation
\begin{equation} \label{eq:dysonr}
  \begin{cases}
    &\paren{i\partial_t-h} \GR(t) - \int_0^t \SR(t-t')
    \GR(t') \, dt' = 0 \\
    &\GR(0) = -i.
  \end{cases}
\end{equation}
for $t \ge 0$, and is taken to be zero for $t<0$.

\subsection{Semi-discretization in the imaginary time variable}

To write the Dyson equation as a system of equations of the form
\eqref{eq:dysongeneral}, we must semi-discretize in the imaginary time
variable $\tau$. In addition to the equation \eqref{eq:dysonm} for
$\GM$, we will obtain $r$ equations from \eqref{eq:dysonlm}, where $r$ is the
number of grid points discretizing $\tau \in [0,\beta]$. It is therefore
important to find as efficient of a discretization of the imaginary
time domain as possible.

Due to the importance of the Matsubara formalism in many-body
calculations, there is a significant literature on the efficient
discretization of $\GM$.
The most common approach is to use a uniform grid on $\tau \in
[0,\beta]$, and to solve \eqref{eq:dysonm} by an FFT-based method. However, since
$\GM$ is not periodic, this discretization suffers from low-order
accuracy, and it struggles to resolve $\GM$ in low temperature
calculations. Methods based on orthogonal polynomial expansions
\cite{boehnke11,gull18,dong20}
and adaptive grids \cite{ku00,ku02,kananenka16} offer an improvement, but remain suboptimal.
Recently, two methods have been proposed which use the specific structure
of imaginary time Green's functions, along with low rank compression
techniques, to obtain highly compact representations: the intermediate representation with sparse sampling
\cite{shinaoka17,chikano18,li20,shinaoka21_2}, and the discrete Lehmann
representation (DLR) \cite{kaye21_3}. Here, we work with the DLR.

Using the DLR, $\GM$ can be expanded in a basis of a small number $r$ of
exponential functions,
\begin{equation} \label{eq:gmdlr}
  \GM(\tau) \approx \sum_{k=1}^r \frac{e^{-\omega_k \tau}}{1+e^{-\beta
  \omega_k}} \wh{g^M_k},
\end{equation}
for carefully chosen DLR frequencies $\omega_k$. The expansion coefficients $\wh{g^M_k}$
can be recovered from values $\gM_j = \GM(\tau_j)$ at a collection of
$r$ DLR imaginary time grid points $\tau_j$. The frequencies
and imaginary time grid points can be obtained using the interpolative
decomposition \cite{cheng05,liberty07}, and they depend only on a dimensionless cutoff
parameter $\Lambda = \beta \omega_{\max}$, with $\omega_{\max}$ a
frequency cutoff, and an accuracy parameter $\epsilon$.
$\Lambda$ can typically be estimated based on physical considerations,
but in practice calculations should be converged with respect to that
parameter. The size of the basis scales as $r = \OO{\log\paren{\Lambda}
\log\paren{1/\epsilon}}$, and is exceptionally small in practice. The
Matsubara component \eqref{eq:dysonm} of the Dyson equation can be solved on the DLR imaginary
time grid. We refer to \cite{kaye21_3} for further details on the
DLR.

Once the DLR expansion of $\GM$ is known, it can be used to compute $\QTV$ by the convolution in
\eqref{eq:qlm}. In particular, let $\sTV_j(t)$ and $\qTV_j(t)$ be the
DLR grid discretizations of $\STV(t,\tau_j)$ and $\QTV(t,\tau_j)$, respectively,
so that $\qTV(t), \sTV(t) \in \CC^r$. Then the discretization
of \eqref{eq:qlm} on the DLR grid is given by
\begin{equation} \label{eq:qlmdisc}
  \qTV(t) = \wb{\GM} \sTV(t),
\end{equation}
for an $r \times r$ matrix $\wb{\GM}$ which can be computed from the
values $\gM_j$. We refer to \cite{kaye21_3} for a detailed
description of this process.
We note that we have assumed $\STV$ and $\QTV$ can also be represented
by the DLR, which is the case for typical problems of physical interest.

As we did for $\sTV$, and $\qTV$, we define $\gTV_j(t)$ as the
discretization of $\GTV(t,\tau_j)$, so that $\gTV(t) \in \CC^R$, and $\sR(t)$ as the discretization of
$\SR(t)$. The nature of the dependence of $\sTV$ and $\sR$ on $\gTV$
follows from that of $\STV$ and $\SR$ on $\GTV$; namely, we have
$\sTV(t) = \sTV(\gTV(t),t)$ and $\sR(t) = \sR(\gTV(t),t)$.

With these definitions, the semi-discretization of \eqref{eq:dysonlm} is
given a coupled system of integro-differential equations,
\begin{equation} \label{eq:dysonlmdisc}
  \begin{cases}
    &(i \partial_t - h) \gTV(t) - \int_0^t \sR(\gTV(t-t'),t-t') \gTV(t') \, dt' =
    \qTV(\gTV(t-t'),t) \\
    &\gTV_j(0) = i \xi \GM(\beta-\tau_j),
  \end{cases}
\end{equation}
for $j = 1,\ldots,r$.
The exponential change of variables $y(t) \equiv e^{i h t} \gTV(t)$ yields a system of equations
of the form \eqref{eq:dysongeneral}, with
\[k(y(t),t) \equiv \sR\paren{e^{-i h t} y(t),t},\] 
\[f(y(t),t) \equiv e^{i h t} \qTV(t) = e^{i h t} \wb{\GM} \sTV(e^{-i h
t} y(t),t),\]
and initial conditions
\[y_j(0) = i \xi \GM(\beta-\tau_j).\]

We note that as for $\GM(\tau)$ in \eqref{eq:gmdlr}, $\GTV(t,\tau)$ can be expanded in the DLR
basis:
\begin{equation} \label{eq:gtvdlr}
  \GTV(t,\tau) \approx \sum_{k=1}^r \frac{e^{-\omega_k \tau}}{1+e^{-\beta
  \omega_k}} \wh{g}_k^\rceil(t).
\end{equation}
The $r$ expansion coefficients $\wh{g}_k^\rceil(t)$ can be obtain from
the $r$ values $\gTV_j(t)$. This expression can be used, for example, to
compute $\GR(t)$ using \eqref{eq:grfromgtv}.

\section{Numerical examples} \label{sec:examples}

As a benchmark and proof of concept we apply the fast history summation
algorithm to two simple but nontrivial systems with different
self-energy expressions: the Bethe graph \cite{QPGFBook,
Georges:1996aa}, for which $\Sigma$ depends linearly on $G$, and the SYK
model \cite{PhysRevLett.70.3339,chowdhury21,Gu:2020vl}, for which the dependence is
cubic.

All numerical experiments were implemented in Fortran, and carried out
on a single CPU core of a laptop with an Intel Xeon E-2176M 2.70GHz
processor. The Fortran library \texttt{libdlr} was used for an
implementation of the DLR \cite{kaye21_libdlr,libdlr}. The FFTW library \cite{frigo05} was used for FFTs.

\subsection{The Bethe graph}

The Bethe graph \cite{QPGFBook, Georges:1996aa} is a simplified model of
a periodic lattice system in which each lattice site is connected to $q$ other
sites.
It is commonly employed in model calculations since the dynamics on the
infinite graph can be described by the simple self-energy
\begin{equation}
  \Sigma = c^2 G.
  \label{eq:bethe}
\end{equation}
Here, the real constant $c = qJ$ is determined by the hopping parameter
$J$, given by the matrix element of the kinetic energy operator coupling
states at neighbouring sites, and the coordination number $q$. The
Green's function is known analytically in this case, so it provides a
straightforward performance test for our time stepping method and
history summation algorithm.

Let us first derive $G^R$ analytically, starting from the Dyson equation
\eqref{eq:dysonr}. We note that, following the convention in
the literature \cite{Negele:1998aa, Stefanucci:2013oq}, we define the
Fourier transform of by $f(\omega) = \int_{-\infty}^\infty e^{i
\omega t} f(t) \, dt$, and the function argument itself indicates
that we are working in the transform domain.
Since by definition $G^R(t) = 0$ and $\Sigma^R(t) = c^2 G^R(t) = 0$ for
$t < 0$, \eqref{eq:dysonr} can be extended to the whole real axis\footnote{
Using the boundary condition $G^R(0) = -i$,
$\Sigma^R(t) = \theta(t) \Sigma^R(t)$, and $G^R(t) = \theta(t) G^R(t)$
gives
\begin{align*}
  i \partial_t G^R(t) &= i \delta(t) G^R(0) + i \theta(t) \partial_t
G^R(t) \\ 
  &= \delta(t) + \theta(t) \left[ hG^R(t) + \int_0^t \Sigma^R(t -
t') G^R(t') \, dt' \right] \\
  &= \delta(t) + hG^R(t) + \int_{-\infty}^\infty
\Sigma^R(t - t') G^R(t') \, dt'.
\end{align*}
}:
\[
\paren{i \partial_t - h} \GR(t) - c^2 \int_{-\infty}^\infty
\GR(t-t') \GR(t') \, dt' = \delta(t).
\]
Applying the Fourier transform yields a quadratic equation for
$\GR(\omega)$,
\[\paren{\omega-h} \GR(\omega) - c^2 (\GR(\omega))^2 = 1,\]
which has the solution
\begin{equation}\label{eq:grtilde1}
  \GR(\omega) = \frac{1}{2c^2} \left[ \omega-h - \sqrt{(\omega-h)^2-4
  c^2} \right].
\end{equation}
Since $G^R(t)$ is zero for $t<0$,
the spectral function $A(\omega) \equiv -\frac{1}{\pi} \text{Im}\,\GR(\omega)$
satisfies the Kramers-Kronig relations,
\[ \GR(\omega) = - i \pi A(\omega)
+
\text{P.V.} \int_{-\infty}^\infty \frac{A(\nu)}{\omega - \nu} \, d\nu =
-i \paren{\theta \ast A}(\omega),
\]
where $\theta(\omega)$ is the Fourier transform of the Heaviside
function.
It follows, in particular, that
\[ G^R(t) = -2 \pi i \theta(t) A(t). \]
From \eqref{eq:grtilde1}, we see that $A(\omega)$
is the well known semi-circular density of states:
\[ A(\omega)
= 
\begin{cases}
  \frac{1}{2\pi c^2} \sqrt{ 4c^2 - (\omega - h)^2} & \text{if } (\omega - h)^2 < 4c^2 \\
  0 & \text{otherwise}.
\end{cases}
\]
Applying the inverse Fourier transform gives
\begin{align*}
  \GR(t)
  &= - i\theta(t) \int_{-\infty}^\infty e^{-i\omega t} A(\omega) \, d\omega \\
  &= -\frac{i \theta(t)}{2 \pi c^2} \int_{h-2c}^{h+2c} e^{-i \omega t} 
    \sqrt{4c^2 - (\omega-h)^2} \, d\omega \\
    &= \frac{\theta(t)}{2 \pi c^2 t} \int_{h-2c}^{h+2c}  
    \frac{(\omega-h) e^{-i \omega t}}{\sqrt{4c^2 - (\omega-h)^2}} \, d\omega
    \\
    &= \theta(t) \frac{e^{-i h t}}{\pi c t}  \int_{-\pi/2}^{\pi/2}
    e^{- i 2c t \sin \theta} \sin \theta \, d\theta,
\end{align*}
where in the third
equality we have integrated by parts, and in the fourth equality we
have used the change of variables $2 c \sin \theta = \omega - h$. Here,
we recognize the integral formula for the Bessel function $J_1$ of the
first kind, and obtain
\begin{equation} \label{eq:grexact}
  \GR(t) = - i \theta(t) e^{-i h t} \frac{J_1(2ct)}{ct}.
\end{equation}

\begin{figure}
  \centering
  \includegraphics[width=0.75\textwidth]{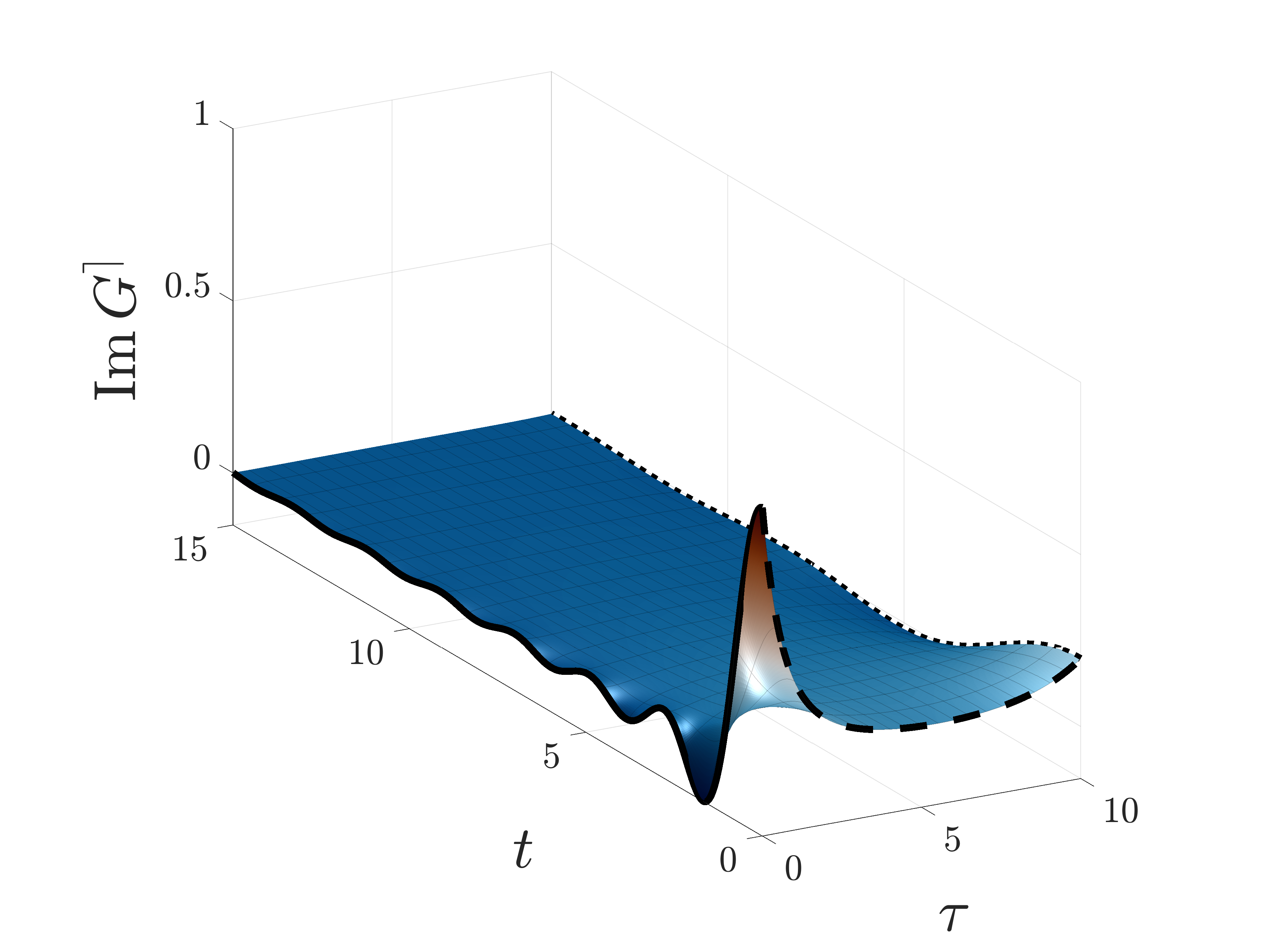}

  \caption{$\Im \GTV(t,\tau)$ for the Bethe lattice self-consistency.
  $\GM(-\tau) = \Im \GTV(0,\tau)$ is indicated by the dashed line,
  $\Im G^<(t) = \Im \GTV(t,0)$ by the solid line, and
  $-\Im G^>(t) = \Im \GTV(t,\beta)$ by the dotted line.}
  \label{fig:gtvbethe}
\end{figure}

For the numerical calculations,
we fix $c = 1$, $h = -1$ and $\beta = 10$, and consider the
fermionic case $\xi = -1$.
Since the initial condition for the real time propagation of $\GTV(t,
\tau)$ in \eqref{eq:dysonlm} is determined by $G^M(\tau)$, we first solve the imaginary time equation of motion \eqref{eq:dysonm}.
The imaginary time axis $\tau \in [0, \beta]$
is discretized using the DLR tolerance $\epsilon = 10^{-15}$,
and a high energy cutoff $\Lambda = 40$, which we have
verified is sufficient to eliminate the discretization error in $\tau$
to high accuracy.
The resulting DLR contains $r = 31$ basis functions
and DLR grid points $\tau_j$.
We obtain $G^M(\tau)$ on the DLR grid using the method described in
\cite{kaye21_3} for the equivalent integral form of \eqref{eq:dysonm}, treating the nonlinearity by a weighted fixed point iteration with pointwise tolerance $10^{-15}$. 

\begin{figure}
  \centering
  \includegraphics[width=0.45\textwidth]{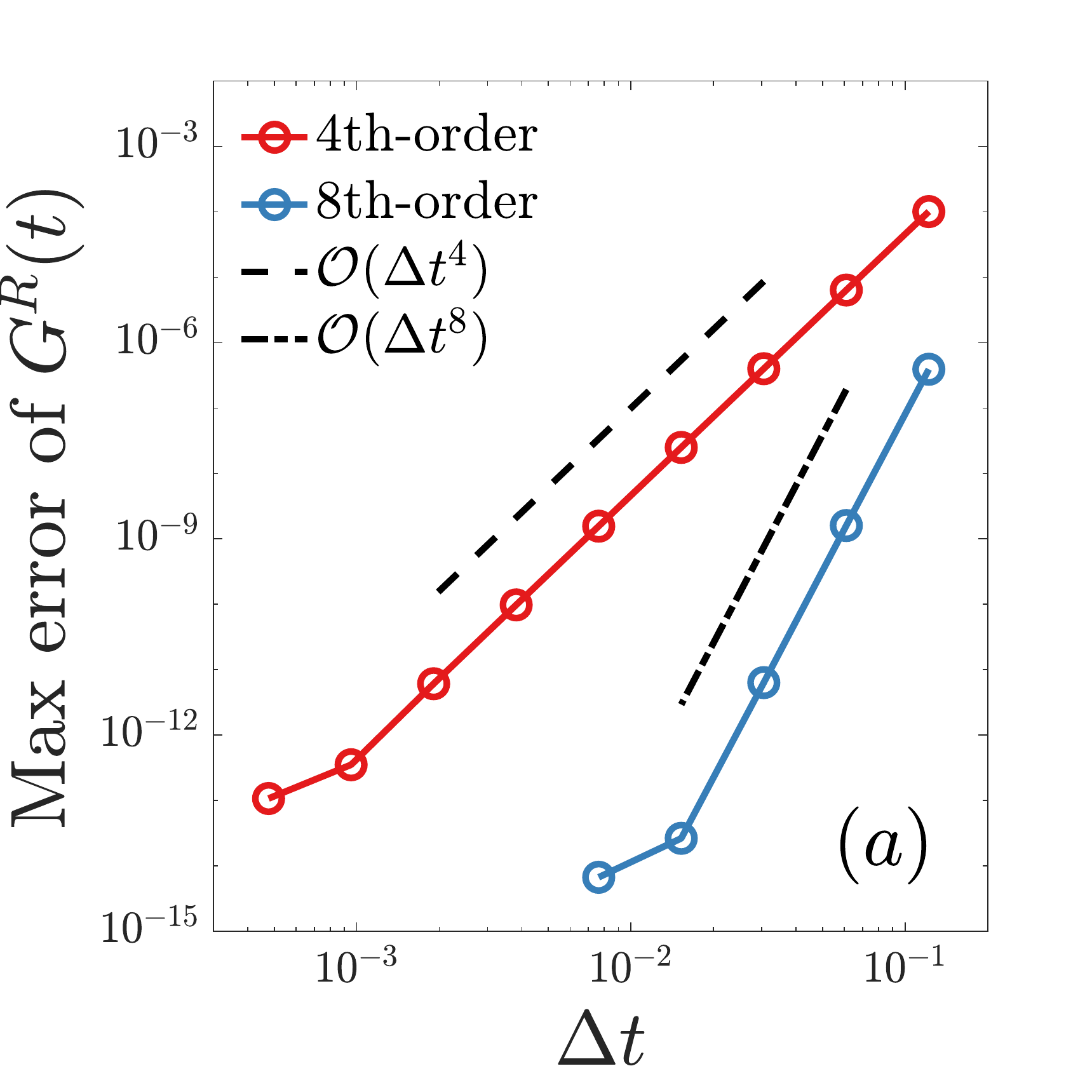}
  \hspace{5ex}
  \includegraphics[width=0.45\textwidth]{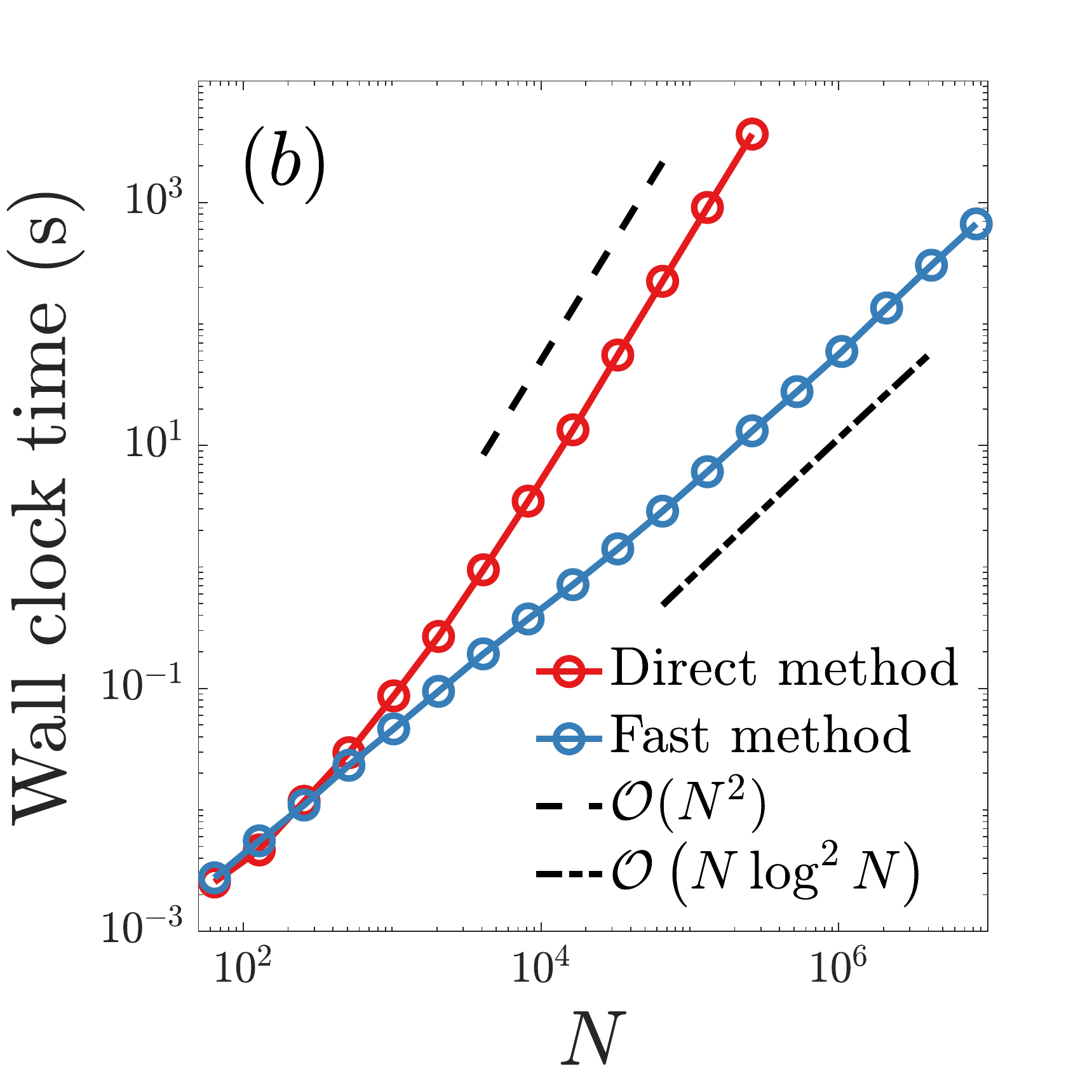}

  \caption{Numerical experiments with the Bethe
  lattice self-consistency $\Sigma = G$. (a) Convergence of the Adams-Moulton method as $\Delta t \to
  0$. (b) Wall clock timings using direct and fast history summation with
  fixed $\Delta t$ and increasing propagation time.}
  \label{fig:bethe}
\end{figure}

We then solve the real time equation of motion \eqref{eq:dysonlm} for $\GTV(t, \tau)$ using the high-order
Adams-Moulton method described in Appendix \ref{app:adams}, and compute
$\GR(t)$ from \eqref{eq:grfromgtv} by evaluating the DLR expansion
\eqref{eq:gtvdlr} of $\GTV$.
We use fixed point iteration with a pointwise tolerance $10^{-15}$ for
the nonlinear solve.
Figure \ref{fig:gtvbethe} shows a plot of $\GTV(t,\tau)$ for $t \in [0,15]$. In
Figure \ref{fig:bethe}a, we show convergence of the fourth and
eighth-order Adams-Moulton methods
with respect to $\Delta t$ by comparing the computed value of $\GR(t)$
with the exact solution \eqref{eq:grexact} for $t \in [0,1000]$. We then
fix $\Delta t = 1/64$, sufficient to achieve near machine precision
accuracy for all calculations using the eighth-order method, and vary the final propagation time
from $t = 1$ to $t = 131072$ by increasing the number $N$ of time steps.
Figure \ref{fig:bethe}b shows wall clock timings using
direct history summation (for some choices of $N$) and our fast history summation
algorithm. We observe the expected $\OO{N \log^2 N}$ scaling for the fast method, and note
that it is superior for $N \geq 256$, indicating a very
small pre-factor in the scaling. For the
longest propagation carried out here, with $N = 8388608$ and $T = 131072$, only
$N \approx 110000$ time steps to $T \approx 1700$ would be possible by the direct
method with the same cost.
We also note that only one fixed point
iteration is required to reach self-consistency for all time steps after
the first $500$, and before that at most two, indicating that the equations are
non-stiff and an Adams predictor-corrector-type method is sufficient in this
case \cite[Sec 5.4.2.]{ascher98}. 

\subsection{The Sachdev-Yitaev-Ke model}

The fermionic Sachdev-Yitaev-Ke (SYK) model demonstrates the challenge of resolving emergent
low energy features in quantum many-body systems. It is used to model
certain features of strange metals, a poorly understood phase
of matter with properties distinct from simple metals, and is connected
to a variety of other phenomena of physical interest
\cite{PhysRevLett.70.3339,chowdhury21,Gu:2020vl,garcia17,wei21}.
The SYK equations of motion in real time are given by \eqref{eq:dysonm} and
\eqref{eq:dysonlm} with $\xi = -1$, and a self-energy consisting of a
single second-order diagram in the interaction $J$, containing a product
of three Green's functions:
\begin{equation} \label{eq:syksigma}
  \begin{aligned}
    \SM(\tau) &= J^2 \paren{\GM(\tau)}^2 \GM(\beta-\tau) \\
    \STV(t,\tau) &= J^2 \paren{\GTV(t,\tau)}^2
    \paren{\GTV(t,\beta-\tau)}^*.
  \end{aligned}
\end{equation}
In the limit of zero temperature $T = 1/\beta \to 0$, the SYK model displays an emergent low energy feature in the spectral
function $A(\omega) =  - \frac{1}{\pi} \Im G^R(\omega)$:
a square root divergence $A(\omega) \sim |\omega|^{-1/2}$
\cite{PhysRevLett.70.3339}. At finite temperature $T$, the divergence is
regularized by thermal fluctuations at frequencies $|\omega| \lesssim T$.

In order to resolve this feature, the Green's function must be
propagated to a large time at low temperature. Our fast history
summation algorithm makes this practical. We demonstrate this by solving
the SYK equations and computing the spectral density at temperatures up
to four orders of magnitude smaller than the coupling strength $J$,
and establish the scaling of the cost required to reach even lower temperatures.

We solve the SYK equations and compute $\GR(t)$ and $A(\omega)$ for $J = 1$, $h = 0$ and
$\beta = 100, 1000, 10000$ using the eighth-order Adams-Moulton method
described in Appendix \ref{app:adams}. We tune
all discretization parameters, and the propagation time, so that the 
spectral function $A(\omega)$ is resolved.

For $\beta = 10000$, we use the DLR parameters $\Lambda = 10^5$ and
$\epsilon = 10^{-10}$, which yields $r=92$ degrees of freedom in $\tau$.
We use a fixed point tolerance of $10^{-14}$ for the nonlinear
iteration. We propagate to a final time $t = 50000$ with $N = 1048576$ time
steps. The computation takes less than five minutes to complete using our fast
method. It would take approximately two days by the direct method.

We note that for all three choices of $\beta$, only one fixed point
iteration is required to reach self-consistency
for all time steps after the first $100$, and before that at most three. 

\begin{figure}
  \centering
  \includegraphics[width=0.45\textwidth]{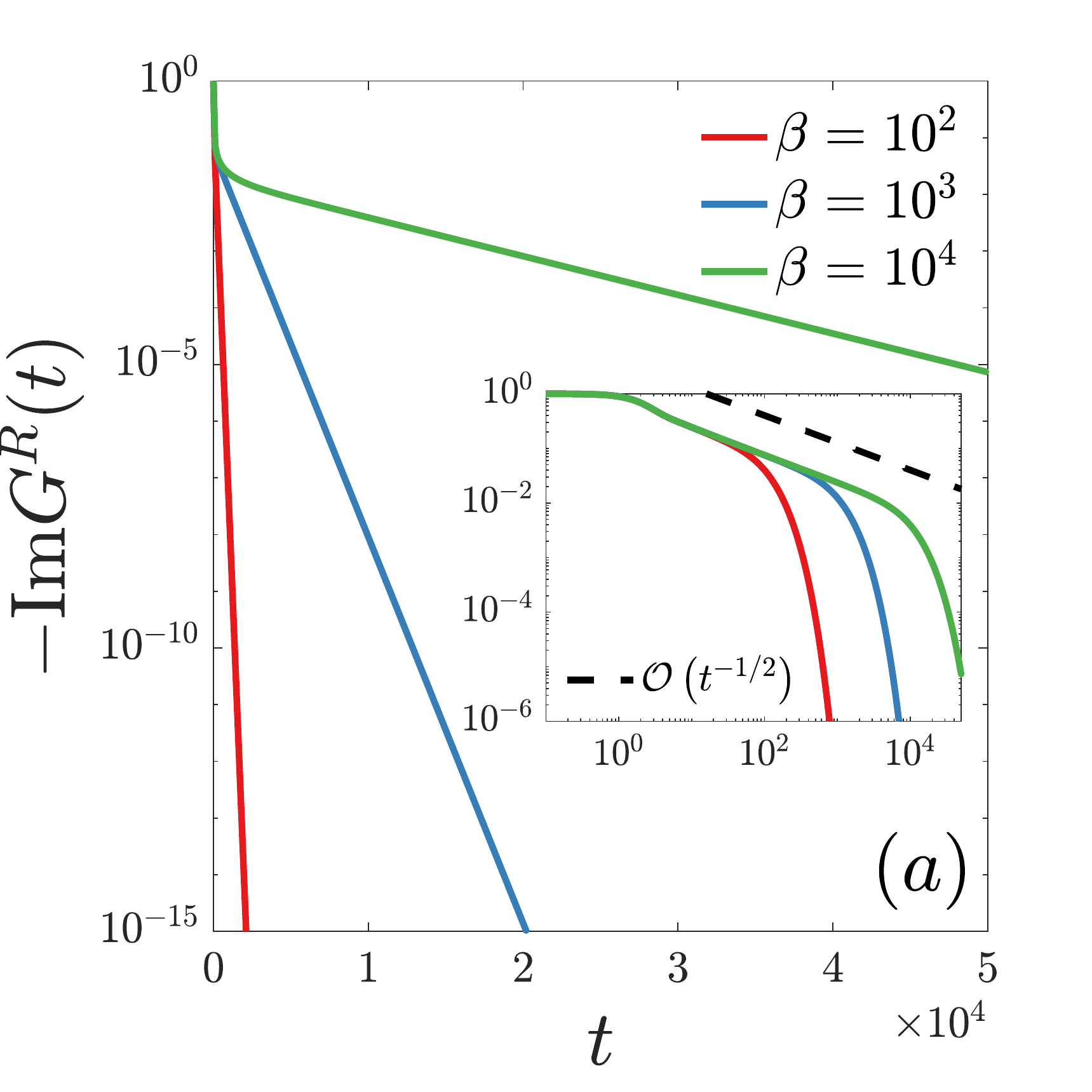}
  \hspace{5ex}
  \includegraphics[width=0.45\textwidth]{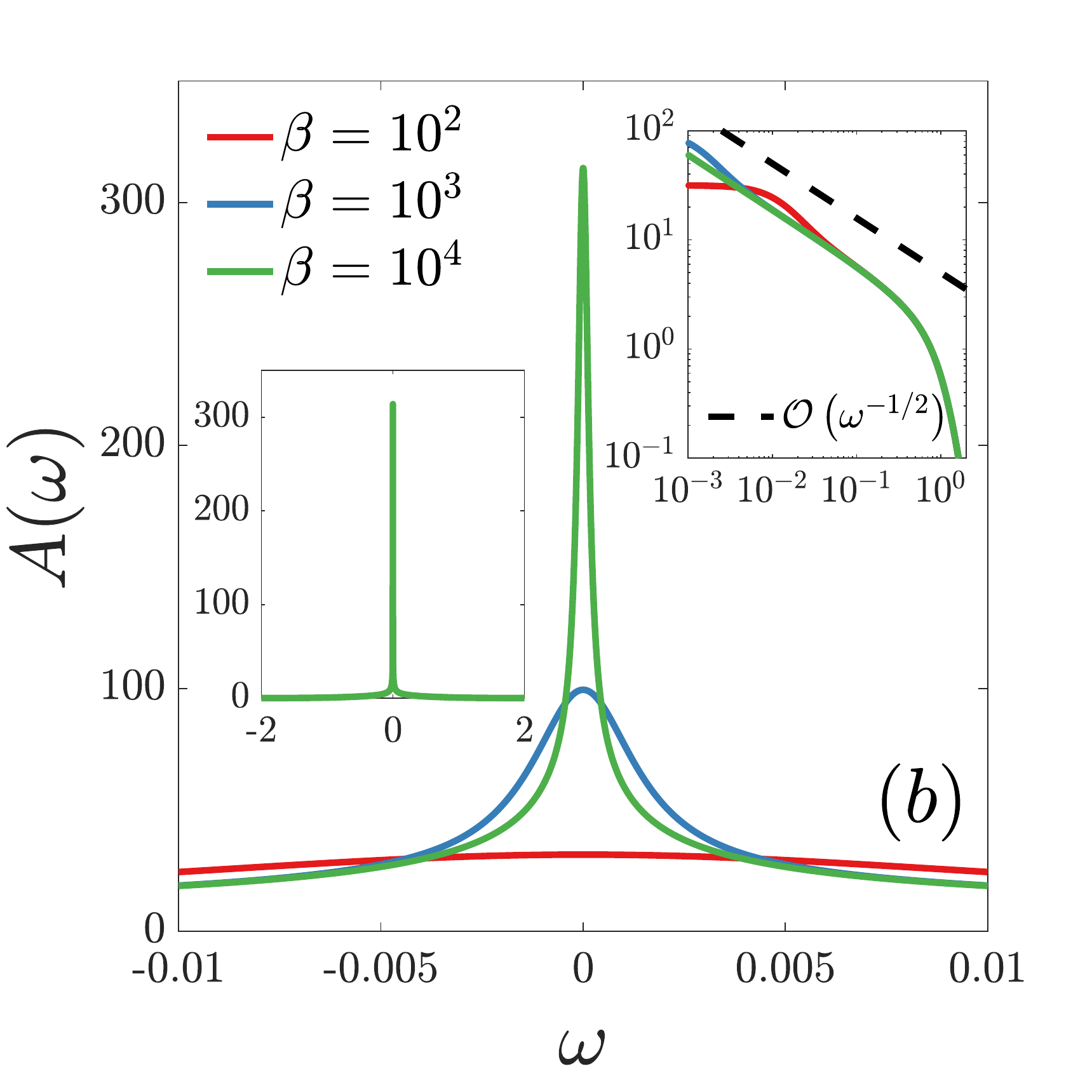}

  \caption{Numerical experiments for the SYK model. (a) Asymptotic
  behavior of the retarded Green's function $G^R(t)$ at different
  temperatures. The inital
  $t^{-1/2}$ decay (inset) transitions into exponential decay at a rate
  proportional to the temperature, requiring large propagation times for
  low temperature. (b) Structure of the spectral density $A(\omega)$
  at different temperatures. The spectral function becomes highly
  localized at low temperatures (left inset). In the $\beta \to \infty$
  limit, it approaches an $|\omega|^{-1/2}$ singularity at the origin (right inset).}
  \label{fig:syk}
\end{figure}

Plots of $\GR$ and $A$ are shown in Figure \ref{fig:syk}. $\GR(t)$ decays
exponentially as $t \to \infty$ at a rate proportional to the
temperature (Figure \ref{fig:syk}a). As a result, for large $\beta$, the
spectral density is highly localized (Figure \ref{fig:syk}b). In the
zero-temperature limit, the spectral density has an
$|\omega|^{-1/2}$ singularity at the origin.  Although for finite
temperature $A(\omega)$ is smooth, we observe it approaching this
singularity as $\beta \to \infty$, with a corresponding
$t^{-1/2}$ decay regime in $G^R(t)$ at intermediate times.

\section{Conclusion} \label{sec:conclusion}

We have presented an efficient solver for a class of nonlinear
convolutional Volterra integro-differential equations
for which the integral kernel depends causally on the
solution. As is typical for Volterra-type equations, the computational bottleneck is the
evaluation of history integrals at each time step.
However, in the present setting, limited knowledge of the kernel at a given time step precludes the
use of the various Volterra history summation techniques which have been
introduced in the literature.
We present an $\OO{N \log^2 N}$ algorithm which properly handles the
structure of the integral operator.

We are motivated to study equations of the form \eqref{eq:dysongeneral} by their appearance
in quantum many-body Green's function methods. In particular, the
equilibrium real time Dyson equation of motion \eqref{eq:dysonlm} for the single particle
Green's function can be reduced to a system of
equations of this form.
By using the imaginary time Green's function, which is relatively easy
to compute, to initialize real time propagation, our approach requires
only local nonlinear iteration. By contrast, the standard method
requires global nonlinear iteration in real time and frequency,
which can lead to slow convergence in regimes of physical interest.

As proofs of concept, we present numerical experiments for the Bethe
graph \cite{QPGFBook, Georges:1996aa} and the Sachdev-Ye-Kitaev model
\cite{PhysRevLett.70.3339, Gu:2020vl}, demonstrating propagation to
extremely large times.
In particular, we are able to resolve the square root divergence in the
SYK model at a temperature $1/\beta$ four orders of magnitude lower than
the coupling $J$, $\beta^{-1}/J = 10^{-4}$, while resolving energy
scales $\omega$ three orders of magnitude smaller than $J$, $\omega / J
\approx 10^{-3}$, in a calculation taking less than five minutes on a
single CPU core of a laptop.
More generally, our approach is applicable to many-body Green's function methods with
low order self-energy expressions, as in low order pertubation theory
\cite{Abrikosov:1975aa, Fetter:2003aa, Negele:1998aa,
Stefanucci:2013oq}, as well as bubble and ladder resummations like
Hedin's GW approximation and beyond \cite{Hedin:1965aa,
10.3389/fchem.2019.00377, PhysRevB.90.115134, PhysRevB.104.085109}.

\section*{Acknowledgements}

We thank Alex Barnett for helpful discussions.
H.U.R.S.\ acknowledges financial support from the ERC synergy grant (854843-FASTCORR).
The Flatiron Institute is a division of the Simons Foundation.

\appendix

\section{Adams-Moulton time stepping with Gregory integration} \label{app:adams}

We consider a general coupled system of equations of the form
\eqref{eq:dysongeneral}, in order to directly treat the case of the
equilibrium Dyson equation described in Section \ref{sec:dyson}:
\begin{equation}\label{eq:dysonsystem}
  i y_j'(t) + \int_0^t
  k_j\paren{y_1(t-t'),\ldots,y_d(t-t'),t-t'} y_j(t') \, dt' =
  f_j(y_1(t),\ldots,y_d(t),t)
\end{equation}
for $j = 1,\ldots,d$.

We discretize \eqref{eq:dysonsystem} by the Adams-Moulton method. 
For notational simplicity, we suppress the dependence of $k$ and $f$ on
$y_1,\ldots,y_d$, and write
\[\wt{f}_j(t) \equiv f_j(t) - \int_0^t k_j(t-t') y_j(t') \, dt'  \]
so that
\[i y_j'(t) = \wt{f}_j(t).\]
Integrating both sides from $t_n$ to $t_{n+1}$ gives
\[y_j^{n+1} = y_j^n - i \int_{t_n}^{t_{n+1}} \wt{f}_j(t) \, dt.\]
The Adams-Moulton method of order $p$ is obtained by replacing the integrand with its
polynomial interpolant at the points $\{t_{n+1-j}\}_{j=0}^{p-1}$, and
computing the resulting integrals analytically. This yields the discretization
\begin{equation} \label{eq:adams1}
  y_j^{n+1} = y_j^n - i \Delta t \sum_{l=0}^{p-1} \muam_l \wt{f}_j^{n+1-l} =
  y_j^n - i \Delta t \sum_{l=0}^{p-1} \muam_l \paren{f_j^{n+1-l}
  - \int_0^{t_{n+1-l}} k_j(t_{n+1-l}-t') y_j(t') \, dt'}.
\end{equation}
Here, $\muam_j$ are the Adams-Moulton weights.
For a procedure
to obtain the weights, along with tabulated weights for
methods of up to eighth-order, we refer the reader to \cite[Chap. 24]{butcher16}.

It remains to discretize the history
integrals to high-order accuracy. Since $y_j(t)$ is known only at the grid points
$t = t_n$, a standard high-order rule, such as a composite
Gauss-Legendre rule, would require interpolation, and is inconsistent
with our fast algorithm without significant modification. Instead we require a high-order rule of the form
\begin{equation} \label{eq:skn}
  \int_0^{t_n} k_j(t_n-t') y_j(t') \, dt' \approx \Delta t \sum_{m=0}^n
  k_j^{n-m} y_j^m + \Delta t \sum_{m=0}^{q-1} \mug_m
\paren{k_j^{n-m} y_j^m + k_j^m y_j^{n-m}} \equiv S_j^n
\end{equation}
for some weights $\mug_m$; that is, a standard equispaced rule with
endpoint corrections. A simple method is the Gregory integration rule. The Gregory rule with
$q$ weights is correct to order $q+1$. Although this approach becomes
unstable for very high orders, we have found it to be effective at least
up to order $8$. More stable endpoint-corrected rules have been
introduced, including modifications of Gregory rules
\cite{fornberg19,fornberg21}, and stable rules of very high order which require modifying endpoint
node locations \cite{kapur97,alpert99}, but we have found the standard Gregory rules to be sufficient for our
purposes. For a procedure to obtain the weights, we refer the reader to
\cite{fornberg21}; the weights are obtained from the Gregory
coefficients, which are tabulated to high order in
\cite{oeisgregn,oeisgregd}.

Inserting \eqref{eq:skn} into \eqref{eq:adams1} and rearranging to place
all unknown quantities on the left hand side, we obtain
\begin{equation} \label{eq:amfulldisc}
  \begin{multlined}
    \paren{1 - i \Delta t^2 \muam_0 \paren{1+\mug_0}}
    y_j^{n+1} - i \muam_0 \Delta t^2 \paren{1+\mug_0} k_j^{n+1} y_j^0 +
    i \muam_0 \Delta t f_j^{n+1} \\ = y_j^n + i \muam_0 \Delta t
    \wb{S}_j^{n+1} - i \Delta t
    \sum_{l=1}^{p-1} \muam_l \paren{f_j^{n+1-l} - S_j^{n+1-l}}
  \end{multlined}
\end{equation}
where we have defined the truncated endpoint-corrected history sum
\begin{equation} \label{eq:trunchist}
  \wb{S}_j^n \equiv S_j^n
  - \Delta t \left(
  k_j^n y_j^0 + k_j^0 y_j^n + \mug_0 \paren{k_j^n y_j^0 + k_j^0 y_j^n}
  \right)
  .
\end{equation}

This is a collection of $d$ coupled nonlinear equations, which are
closed by prescribing the functional dependence of $k_j^n$ and $f_j^n$
on $y_1^n,\ldots,y_d^n$.
Once these are solved, $S_j^{n+1}$ must be
computed by updating the value $\wb{S}_j^{n+1}$
using \eqref{eq:trunchist}. In any given time
step, the previous $p-1$ computed values of $f_j^n$ and $S_j^n$ must be
stored, along with the full history of $y_j$, for $j=1,\ldots,d$.

We can split the endpoint corrected history sums $S_j^n$ into two terms, $S_j^n =
\Delta t \paren{s_j^n + c_j^n}$, with
\[s_j^n \equiv \sum_{m=0}^n k_j^{n-m} y_j^m,\]
a standard equispaced history sum, and
\[c_j^n \equiv \sum_{m=0}^{q-1} \mug_m
\paren{k_j^{n-m} y_j^m + k_j^m y_j^{n-m}},\]
a local correction. Thus, as for the simple discretization given in
Section \ref{sec:fastalg}, we recognize the computation of the sums
$s_j^n$ as the quadratic scaling bottleneck of the time stepping
procedure, and we can apply the same fast algorithm.

We do not suggest a specific nonlinear iteration procedure to solve
\eqref{eq:amfulldisc}
-- in our numerical examples, we use a simple fixed point iteration -- but
propose using an Adams-Bashforth method to obtain an initial guess. The
Adams-Bashforth method is derived in a similar manner to the
Adams-Moulton method, but it is explicit. The $p$th order
Adams-Bashforth discretization is given by 
\begin{equation} \label{eq:ab}
  \begin{multlined}
    y_j^{n+1} = y_j^n - i \Delta t \sum_{l=1}^p \muab_l \paren{f_j^{n+1-l} -
    S_j^{n+1-l}}.
  \end{multlined}
\end{equation}
We again refer the reader to \cite[Chap. 24]{butcher16} for a method to
compute the Adams-Bashforth weights $\muab$, along with tabulated
weights for the methods up to eighth-order.
In practice, using this guess can significantly reduce the number of
nonlinear iterations required at each time step.

To obtain a method of order $p$, it suffices to set $q = p-1$. Then
\eqref{eq:amfulldisc} is only applicable when $n \geq p-2$, and
\eqref{eq:ab} when $n \geq p-1$. An
alternative method is then required to obtain $y_j^n$, $f_j^n$, and $S_j^n$
for $n = 0,\ldots,p-1$. A convenient approach is to use Richardson
extrapolation on the second-order version of the method -- the implicit
trapezoidal rule -- which does not
require any initialization. The method, which requires $p$ to be even,
proceeds as follows:
\begin{enumerate}
  \item Use the second-order method, \eqref{eq:amfulldisc}, with $p =
    2$, until a final time $t = (p-1) \Delta t$, taking time steps of size
$\Delta t$, $\Delta t / 2$, $\Delta t / 4$, \ldots, $\Delta t /
    2^{p/2-1}$. Record the computed values of $y_j^n$, $f_j^n$, and
    $S_j^n$ for each time step choice.
  \item Apply Richardson extrapolation to these values to
obtain approximations of $y_j(t)$, $f_j(t)$, and $S_j(t)$ at $t = \Delta t, 2 \Delta t, \ldots,
(p-1) \Delta t$ which are accurate to order $p$.
\end{enumerate}
Here, we have used that the truncation error of the implicit trapezoidal
rule contains only even-order terms. For further details of the
Richardson extrapolation procedure, we refer the reader to \cite[Sec. 3.4.6]{dahlquist08}.

\section{Fast application of square, triangular, and
parallelogram-shaped Toeplitz blocks} \label{app:toeplitz}

We first review the standard fast algorithm to apply an $n \times n$ Toeplitz matrix $A$
in $\OO{n \log n}$ operations \cite[Sec. 4.7.7]{golub96}. Let $D$ be the $n
\times n$
discrete Fourier transform matrix, $D_{jk} = \omega_n^{jk}$ with
$\omega_n = e^{-2 \pi i/n}$, and let $C$ be an $n \times n$ circulant matrix with first
column $c$. Then $D$ diagonalizes $C$:
\[C = D^{-1} \Phi D.\]
Here $\Phi$ is the diagonal matrix with entries $Dc$. Since
$D$ and $D^{-1}$ can be applied to a vector in $\OO{n \log n}$
operations using
the FFT, this gives an $\OO{n \log n}$ algorithm to compute a matrix-vector product
$b = Cx$:
\begin{algorithm}[H]
  \begin{algorithmic}[1]
    \State Compute $Dc$ by FFT
    \State Compute $Dx$ by FFT
    \State Compute the entrywise product $Dc \odot Dx$
    \State Compute $b = D^{-1} \paren{Dc \odot Dx}$ by inverse FFT
  \end{algorithmic}
  \caption{Fast matrix-vector product $b = Cx$ for a circulant matrix $C$
  with first column $c$}
  \label{alg:fastcirc}
\end{algorithm}

To obtain an $\OO{n \log n}$ algorithm to compute the product $b = Ax$
with an $n \times n$ Toeplitz matrix $A$, we simply embed $A$ in a
$2n \times 2n$ circulant matrix $C$, apply $C$ to a zero-padded vector
$x$ using Algorithm \ref{alg:fastcirc}, and truncate the result. An example for $n = 4$ illustrates the
method:
\[
  \left(
  \begin{array}{cccc|cccc}
    \bf{d} & \bf{c} & \bf{b} & \bf{a} & 0 & g & f & e \\
    \bf{e} & \bf{d} & \bf{c} & \bf{b} & a & 0 & g & f \\
    \bf{f} & \bf{e} & \bf{d} & \bf{c} & b & a & 0 & g \\
    \bf{g} & \bf{f} & \bf{e} & \bf{d} & c & b & a & 0 \\
    \hline 
    0 & g & f & e & d & c & b & a \\
    a & 0 & g & f & e & d & c & b \\
    b & a & 0 & g & f & e & d & c \\
    c & b & a & 0 & g & f & e & d \\
  \end{array}
  \right)
  \left(
  \begin{array}{c}
    \bf{x_1} \\
    \bf{x_2} \\
    \bf{x_3} \\
    \bf{x_4} \\
    \hline 
    0 \\
    0 \\
    0 \\
    0 \\
  \end{array}
  \right)
  =
  \left(
  \begin{array}{c}
    \bf{b_1} \\
    \bf{b_2} \\
    \bf{b_3} \\
    \bf{b_4} \\
    \hline 
    - \\
    - \\
    - \\
    - \\
  \end{array}
  \right).
\]
Here, the top-left block, with bolded entries, is the original Toeplitz matrix
$A$, which has been extended to a circulant matrix. The vector $x$ is zero-padded
as shown. Then, the result is truncated; here, the desired result is the
vector $b = (b_1\,\,b_2\,\,b_3\,\,b_4)^T$, and the dashes indicate entries
that are ignored.

A triangular Toeplitz matrix is simply a special case of a square
Toeplitz matrix; in the above example, we would take $e=f=g=0$.
Therefore the algorithm for this case is the same.

We could also treat the case of a Toeplitz matrix whose non-zero entries
form a parallelogram as simply another special case. However, this
requires forming a circulant matrix of approximately twice the size
necessary. We illustrate a more efficient method using the example of a $4 \times 7$
parallelogram:
\[
  \left(
  \begin{array}{cccccccc}
    \bf{d} & \bf{c} & \bf{b} & \bf{a} & 0 & 0 & 0 \\
    0 & \bf{d} & \bf{c} & \bf{b} & \bf{a} & 0 & 0 \\
    0 & 0 & \bf{d} & \bf{c} & \bf{b} & \bf{a} & 0 \\
    0 & 0 & 0 & \bf{d} & \bf{c} & \bf{b} & \bf{a} \\
    \hline 
    a & 0 & 0 & 0 & d & c & b \\
    b & a & 0 & 0 & 0 & d & c \\
    c & b & a & 0 & 0 & 0 & d \\
  \end{array}
  \right)
  \left(
  \begin{array}{c}
    \bf{x_1} \\
    \bf{x_2} \\
    \bf{x_3} \\
    \bf{x_4} \\
    \bf{x_5} \\
    \bf{x_6} \\
    \bf{x_7} \\
  \end{array}
  \right)
  =
  \left(
  \begin{array}{c}
    \bf{b_1} \\
    \bf{b_2} \\
    \bf{b_3} \\
    \bf{b_4} \\
    \hline 
    - \\
    - \\
    - \\
  \end{array}
  \right).
\]
Here, the parallelogram has been embedded in the top half of a circulant
matrix. No zero-padding of the input vector is necessary, but the result
is still truncated to obtain $b$.

\bibliographystyle{ieeetr}
{\footnotesize \bibliography{dysonrt}}

\end{document}